\documentclass[12pt,a4paper,reqno]{amsart}   
 
\usepackage{footmisc}
\usepackage{footnote}  
\usepackage{geometry}
\usepackage{mathtools} 
\usepackage{url}
\usepackage[hidelinks,breaklinks, urlcolor=black]{hyperref}
 
\usepackage[headings]{fullpage}
\usepackage{datetime}   
\usepackage[english]{babel} 
\usepackage{enumerate}  

\usepackage[full]{textcomp} 
\usepackage{newtxtext} 
\usepackage{cabin} 
\usepackage{zlmtt}
\usepackage{mathptmx} 
\usepackage[cal=boondoxo]{mathalfa} 

\usepackage{microtype} 

\newcommand{\Q}{\mathbb{Q}}
\newcommand{\N}{\mathbb{N}}

\newcommand{\Z}{\mathbb{Z}}
\newcommand{\G}{\mathfrak{G}}

\newcommand{\Odi}[1]{\Odip{}{#1}}

\newcommand{\Odip}[2]{\mathcal{O}_{#1}\left(#2\right)}

\allowdisplaybreaks

\renewcommand{\qedsymbol}{$\square$}

\newtheoremstyle{sltheorems}
{10pt}
{6pt}
{\slshape}
{}
{\bfseries}
{.}
{.5em}
{\thmname{#1}\thmnumber{ #2}\thmnote{ (#3)}}

\theoremstyle{sltheorems}

\newtheoremstyle{remark}
{10pt}
{6pt}
{\rm} 
{}
{\bfseries}
{.}
{.5em}
{\thmname{#1}\thmnumber{ #2}\thmnote{ (#3)}}
 \theoremstyle{remark}

\usepackage{caption} 
\allowdisplaybreaks

\begin{document} 

\title[Kummer ratio for  $\Q(\zeta_q)$]{On the computation of the Kummer ratio \\ of the class number for prime cyclotomic fields$^\star$}

 \author[A.~Languasco]{Alessandro Languasco$\dag$}   
\thanks{\mbox{}\hskip-0.375truecm$\dag$: \emph{corresponding author}; Universit\`a di Padova,
 Dipartimento di Matematica,  ``Tullio Levi-Civita'',
Via Trieste 63, 35121 Padova, Italy.
\emph{email}: alessandro.languasco@unipd.it}

\author[P.~Moree] {Pieter Moree$^*$}     

 \author[S.~Saad Eddin]{Sumaia Saad Eddin$\ddag$}       
\thanks{$\ddag$: Institute of Financial Mathematics and Applied Number Theory,
Johannes Kepler University,
Altenbergerstrasse 69, 4040 Linz, Austria. \emph{email}: Sumaia.Saad\_Eddin@jku.at}
  
 \author[A.~Sedunova]{Alisa Sedunova$^*$}   
  \thanks{$^*$: Max-Planck-Institut f\"ur Mathematik,
Vivatsgasse 7, D-53111 Bonn, Germany. \emph{email}: moree@mpim-bonn.mpg.de, alisa.sedunova@phystech.edu}

 \date{\today, \currenttime} 

\begin{abstract}  
Let  $\zeta_q$ 
be a primitive $q^{\text{th}}$ root of unity with $q$ 
an arbitrary odd prime. 
The ratio 
of Kummer's first factor of the class number of the 
cyclotomic number field $\Q(\zeta_q)$ and its 
expected order of magnitude (a simple function of $q$)
is called the Kummer ratio and denoted by $r(q)$.
It is known that typically $r(q)$ is close to 1, but 
nevertheless it is
believed that it is unbounded, but only large on a very thin sequence of primes $q$.

We propose an algorithm to compute $r(q)$ 
requiring the evaluation of $\Odi{q\log q}$ products 
and $\Odi{q}$ logarithms.
Using it we obtain a new record maximum for  $r(q)$, namely $r(6766811) =1.709379\dotsc$ 
(the old record being $r(5231)=1.556562\dotsc$).
The program used and the results described here, are collected 
at the following address \url{http://www.math.unipd.it/~languasc/rq-comput.html}.
 
\smallskip
\textbf{$^\star$: This is a (preliminary) report about the computational part 
of a joint project with Pieter Moree, Sumaia Saad Eddin and Alisa Sedunova}.
\end{abstract} 
\maketitle
\makeatletter
\def\subsubsection{\@startsection{subsubsection}{3}%
  \z@{.3\linespacing\@plus.5\linespacing}{-.5em}%
  {\normalfont\bfseries}}
\makeatother
\section{introduction}
Let  $\zeta_q$ 
be a primitive $q^{\text{th}}$ root of unity with $q$ 
an arbitrary odd prime.
Put
\[  
G(q) := 2q \Bigl(\frac{q}{4\pi^2}\Bigr)^{\frac{q-1}{4}}.
\]
The function $G(q)$ is of super-exponential 
growth in $q$; for example
$G(439) \approx 10^{117}$, $G(3331) \approx 10^{1607}$,
$G(9689) \approx 10^{5792}$
and
$G(2918643191) \approx 10^{5741570411}$.

Kummer \cite{Kummer1851} proved that $h_1(q)$, the first factor of the class number of
the cyclotomic number field $\Q(\zeta_q)$, is a positive integer and conjectured  that 
$h_1(q) \sim G(q)$ as $q \to +\infty$. 
We define the \emph{Kummer ratio} 
as
\begin{equation}
\label{rq-def}
r(q) 
:=
\frac{h_1(q)}{G(q)}.
\end{equation}

\section{First method: using the digamma function}
\label{Hasse-method} 
Recall that $q$ is an odd prime
and let $\chi$ be a primitive odd Dirichlet character mod $q$. 
Using  Hasse's theorem \cite{Hasse1952}
we have that 
\[
 h_1(q) = G(q)  \prod_{\chi\,\textrm{odd}} L(1,\chi)
\]
and, by \eqref{rq-def},  it follows that
\begin{equation}
\label{basic-formula-L}
r(q)  
=
\prod_{\chi\, \textrm{odd}} L(1,\chi).
\end{equation}
Recalling  eq.~(3.1) of \cite{FordLM2014}, \emph{i.e.},
\begin{equation}
\label{L-psi}
L(1,\chi)
=
-
\frac{1}{q}
\sum_{a=1}^{q-1} \chi(a)\ 
\psi \bigl(\frac{a}{q}\bigr), 
\end{equation}
 where $\psi(x)=(\Gamma^\prime/\Gamma)(x)$ is the \emph{digamma} function,
inserting \eqref{L-psi} into \eqref{basic-formula-L} we can also write
\begin{equation}
\label{psi-formula}
r(q) 
=   
\Bigl(\frac{-1}{q}\Bigr)^{\frac{q-1}{2}}
\prod_{\chi\, \textrm{odd}}  
\sum_{a=1}^{q-1} \chi(a)\ 
\psi \bigl(\frac{a}{q}\bigr).
\end{equation}
%
%
Computationally it is more convenient to work with
$\log r(q)$ rather than $r(q),$ which leads to
\begin{equation} 
\label{Hasse-method-formula}
\log r(q)
= 
\frac{q-1}{2} \Bigl(i \pi -  \log q \Bigr) 
+
\sum_{\chi\, \textrm{odd}} 
\log \ \Bigl( 
\sum_{a=1}^{q-1} \chi(a)\ 
\psi \bigl(\frac{a}{q}\bigr)
\Bigr),
\end{equation}
in which  the last logarithm is a complex one. 
It is clear   that  the  sum over odd Dirichlet characters in \eqref{Hasse-method-formula} has an imaginary 
part equal to $- \pi(q-1)/2$; hence  
\begin{equation}
\label{Hasse-method-formula2}
\log r(q)
= 
 -  \frac{(q-1)}{2} \log q 
+
\sum_{\chi\, \textrm{odd}}  
\log \ \Bigl\vert
\sum_{a=1}^{q-1} \chi(a)\ 
\psi \bigl(\frac{a}{q}\bigr)
\Bigr\vert .
\end{equation}

Since in this formula only \emph{odd} Dirichlet characters appear, we can embed a \emph{decimation in frequency 
strategy}\footnote{We use here this nomenclature since it
is standard in the literature on the Fast Fourier Transform, but clearly it can be translated into number theoretic
language using suitable properties of Dirichlet characters.} 
in the Fast Fourier Transform (FFT) algorithm  to perform the sum over $a$, see section \ref{Kummerq-comput}, thus  replacing the  digamma function  with the cotangent one.

Recalling that FFT on a vector of length $N$ costs 
$\Odi{N\log N}$ arithmetic operations, we get  
 that  computing $r(q)$ via \eqref{Hasse-method-formula2}
has a computational cost of $\Odi{q\log q}$ products plus the cost of computing 
$(q-1)/2$ values of  the cotangent and  logarithm functions. This is good but, at least from a theoretical point of view,
we can do slightly better, see the next  section,
even if the performances in the  practical computations of Table \ref{table2} below  are very similar.

\section{Second method: using a generalized Bernoulli number}
\label{chi-Bernoulli-method} 

Let $\chi$ be a primitive odd Dirichlet character mod $q$ with $q$ 
an odd prime.
We define the first $\chi$-Bernoulli number 
$B_{1,\chi}$ (see Proposition 9.5.12 of Cohen \cite{Cohen2007}) as 
\begin{equation}
\label{chi-Bernoulli-def}
B_{1,\chi}
:=
\frac{1}{q} \sum_{a=1}^{q-1}  a \chi(a) .
\end{equation}
By eq.~(2.1) of Shokrollahi \cite{Shokrollahi1999} 
we have 
\begin{equation}
\label{first-factor-Bernoulli}
h_1(q) = 2q \prod_{\chi\, \textrm{odd}} \Bigl(- \frac{B_{1,\chi}}{2} \Bigr).
\end{equation}
Inserting \eqref{chi-Bernoulli-def} and \eqref{first-factor-Bernoulli} into   \eqref{rq-def},  we obtain
\begin{equation}
\label{Bernoulli-formula}
r(q) 
=
  \Bigl(\frac{q}{4\pi^2}\Bigr)^{-{\frac{q-1}{4}} }\prod_{\chi\, \textrm{odd}} \Bigl(- \frac{B_{1,\chi}}{2} \Bigr)
  =
  \Bigl(-\frac{\pi}{q^\frac32}\Bigr)^{\frac{q-1}{2}}  
  \prod_{\chi\, \textrm{odd}}    \sum_{a=1}^{q-1}  a \chi(a),
\end{equation}
which leads to
\begin{equation}
\label{chi-Bernoulli-method-formula}
\log r(q)
=  
\frac{q-1}{2} \Bigl( \log \pi - \frac{3}{2} \log q + i \pi\Bigr) 
+
\sum_{\chi\, \textrm{odd}} 
\log\Bigl( \sum_{a=1}^{q-1}  a \chi(a) \Bigr),
\end{equation}
in which  the last logarithm is a complex one. Moreover it is clear  that the
 sum over the odd Dirichlet characters in \eqref{chi-Bernoulli-method-formula} has an imaginary 
part equal to $- \pi(q-1)/2$; hence 
\begin{equation}
\label{chi-Bernoulli-method-formula2}
\log r(q)
=  
\frac{q-1}{2} \Bigl( \log \pi - \frac{3}{2} \log q \Bigr) 
+ 
\sum_{\chi\, \textrm{odd}} 
\log \ \Bigl\vert
\sum_{a=1}^{q-1}  a \chi(a) 
\Bigl\vert.
\end{equation}
This formula for 
$\log r(q),$ unlike  \eqref{Hasse-method-formula2}, does not require the computation of
values of a special function: it is enough to use the sequence $a=1,\dotsc, q-1$.  
Since only the sum over odd Dirichlet characters is needed, we can embed a \emph{decimation in frequency 
strategy} in the Fast Fourier Transform (FFT) algorithm  to perform the sum over $a$, see section \ref{Kummerq-comput}.
As in the previous section,
it is easy to see that computing $r(q)$ via \eqref{chi-Bernoulli-method-formula2}
has a computational cost of $\Odi{q\log q}$ 
arithmetic operations plus the cost of computing 
$(q-1)/2$ values of  the  logarithm function and products:  so far, this is
the fastest known algorithm  to compute $r(q)$.

This way we were able to get a new record maximal value for  $r(q)$, namely
\[
r(6766811) =1.709379042\dotsc,
\]
see Table \ref{table3}; such a  result was also double-checked using the  method  of section \ref{Hasse-method}. 
The previously known record $r(5231)=1.556562\dotsc$ 
is due to Shokrollahi \cite{Shokrollahi1999}.
We will see more on these computations in section \ref{Kummerq-comput}.

\section{Comparing methods, results and running times} 
\label{Kummerq-comput}

First of all we notice that PARI/Gp, v.~2.11.2, has the ability to  generate
 Dirichlet $L$-functions (and  many other $L$-functions)
and hence, using \eqref{basic-formula-L}, 
the value of $\log r(q)$ $= \sum_{\chi\, \textrm{odd}}$ $\log  \vert L(1,\chi)\vert$
can be  obtained 
with few instructions of the gp scripting language.
This computation has a linear cost in the number of calls
of the {\tt lfun} function of PARI/Gp and it is, at least
on our Dell Optiplex desktop machine, slower than 
both the approaches we are about to describe below.
 
 The other approaches we can use  to compute $\log r(q)$  are the following:
\begin{enumerate}[a)]
\item 
\label{Hasse-approach}
use formula \eqref{Hasse-method-formula} and   the $\psi$-values;  
\item 
\label{Bernoulli-approach}
use  formula \eqref{chi-Bernoulli-method-formula} and   the first  $\chi$-Bernoulli number.
\end{enumerate}

 This way we can double check the computation we will perform.   
 In both   \eqref{Hasse-method-formula}  and \eqref{chi-Bernoulli-method-formula} we remark that,
since $q$ is prime, it is enough to 
determine a primitive root $g$ of $q$, which leads to
the Dirichlet character $\chi_1$ mod $q$ uniquely
determined by 
 $\chi_1(g) = e^{2\pi i/(q-1)}.$ The 
 set of non-trivial characters
 mod $q$ is then $\{\chi_1^j \colon j=1,2,\dotsc,q-2\}$.
 Hence, if, for every $k\in \{0,\dotsc,q-2\}$, we denote $g^k\equiv a_k\in\{1,\dotsc,q-1\}$,
 every summation in  \eqref{Hasse-method-formula}  and \eqref{chi-Bernoulli-method-formula}
 is  of the type $\sum_{k=0}^{q-2}  e^{2\pi i j k /(q-1)} f(a_k/q)$, where $j\in\{1,\dotsc,q-2\}$
 is odd  and $f$ is a suitable function. 
As a consequence, such quantities are
 the Discrete Fourier Transform (DFT)  of the sequence $\{ f(a_k/q)\colon k=0,\dotsc,q-2\}$. 
This observation is due to Rader \cite{Rader1968} and it was used in \cite{FordLM2014} and in \cite{Languasco2019}
 to speed up the computation of similar quantities via the use of Fast Fourier Transform dedicated
 software libraries. 
 
 In our case we can also use the \emph{decimation in frequency} strategy: following the line in section
 4.1 of \cite{Languasco2019}, letting $e(x):=\exp(2\pi i x)$, $m=(q-1)/2$,    
 for every $j=0,\dotsc, q-2$, $j=2t+\ell$, $\ell\in\{0,1\}$ and $t\in \Z$, we have that  
\begin{align}
\notag
 \sum_{k=0}^{q-2}   e\Bigl(\frac{ j k}{q-1}\Bigr)  f \Bigl(\frac{a_k}{q}\Bigr)
 &=
 \sum_{k=0}^{m-1}  
 e\Bigl(\frac{ t k}{m}\Bigr)  
    e\Bigl(\frac{\ell k}{q-1}\Bigr)    
 \Bigl(
 f\Bigl(\frac{a_k}{q}\Bigr) 
 +
 (-1)^{\ell} 
 f \Bigl(\frac{a_{k+m}}{q}\Bigr)
 \Bigr)
 \\&
  \label{DIF} 
 =
 \begin{cases}
 \sum\limits_{k=0}^{m-1}     e\bigl(\frac{t k}{m}\bigr) b_k 
 & \textrm{if} \ \ell =0;\\
  \sum\limits_{k=0}^{m-1}     e\bigl(\frac{t k}{m}\bigr)  c_k  
 & \textrm{if} \ \ell =1,\\
 \end{cases}
\end{align}
where $t=0,\dotsc, m-1$, 
\[
b_k :=
  f\Bigl(\frac{a_k}{q}\Bigr) +  f \Bigl(\frac{a_{k+m}}{q}\Bigr)   
\quad
\textrm{and}
\quad
c_k :=  
   e\Bigl(\frac{k}{q-1}\Bigr)   
 \Bigl(  f\Bigl(\frac{a_k}{q}\Bigr) -  f \Bigl(\frac{a_{k+m}}{q}\Bigr)  \Bigr).
\]
Since we just need the sum  over the odd Dirichlet characters for $f(x)=x$ of $f(x)=\psi(x)$,
instead of computing an FFT transform of length $q-1$
we can evaluate  an FFT of half a length, applied on a suitably
modified sequence according to \eqref{DIF}.
Clearly this leads to a gain in speed and a reduction in memory occupation
in running the actual computer program. 
In case $f(x)= \psi(x)$ we can simplify the expression
$e(  k / (q-1))
 \bigl(\psi(a_k/q) - \psi(a_{k+m}/q) \bigr)$ for $c_k,$ where $m=(q-1)/2$ and $k=0,\dotsc, m-1$,
in the following way. Recalling that $\langle g \rangle = \Z^*_q$, $a_k \equiv g^k \bmod q$ and   $g^m \equiv q-1 \bmod{q}$,
 we can write 
\[ 
\psi \Bigl(\frac{a_{k+m}}{q} \Bigr)  = \psi\Bigl(\frac{q-a_{k}}{q}\Bigr)
=
\psi\Bigl(1-\frac{a_{k}}{q}\Bigr)
\] 
and hence, using the well-known \emph{reflection formula} $\psi(1-x) - \psi(x)  = \pi \cot (\pi x)$,  
we obtain
\begin{align*}
\psi\Bigl(\frac{a_{k}}{q}\Bigr)
-
\psi\Bigl(1-\frac{a_{k}}{q}\Bigr)
&=
-\pi \cot \Bigl(\frac{\pi a_{k}}{q}\Bigr),
\end{align*}
for every $k=0,\dotsc, m-1$. Inserting the last relation in the definition of $c_k$  in \eqref{DIF}
we can replace in the actual computation  the  digamma function  with the cotangent one.
The case $f(x)=x$ is easier; using again $\langle g \rangle = \Z^*_q$, $a_k \equiv g^k \bmod q$ and   $g^m \equiv q-1 \bmod{q}$,
 we can write that $ a_{k+m}  \equiv  q-a_{k} \bmod{q}$;  hence
\[
 a_k  -    a_{k+m}  
 = 
  a_k -(q-a_{k})  
=
2a_k  -q
\]
so that in this case we obtain $c_k=   e(k / (q-1))(2a_k/q -1)$ for every $k=0,\dotsc m-1$, $m=(q-1)/2$.

\subsection{Computations with trivial summing over $a$ (slower, more decimal digits available).}
 
 Unfortunately in \texttt{libpari} the FFT-functions work only if $q=2^\ell+1$, for some $\ell\in \N$. 
So  we had to trivially perform these summations and  hence, in practice,  
this part is the most time consuming one in both approach  
\ref{Hasse-approach}) and \ref{Bernoulli-approach}), as the cost is
quadratic in $q$.

Being aware of such limitations, 
we performed the computation of $r(q)$ with  these three approaches for every $q$ prime, $3\le q\le 1000$, 
on a Dell OptiPlex-3050, equipped with an Intel i5-7500 processor, 3.40GHz, 
16 GB of RAM and running Ubuntu 18.04.2, using a   
precision of $30$ decimal digits, see Table \ref{table1}. 
The results  coincide up the desired 
precision.
We also computed the values of $r(q)$ for   
$q=1451$, $2741$, $3331$, $4349$, $4391$, $5231$, $6101$, $6379$, $7219$, $8209$, $9049$, $9689$
and these are given in Table \ref{table2}.
These numbers were chosen to  extend the number of available decimals for the known data (see Fung-Granville-Williams \cite{FungGW1992} and Shokrollahi \cite{Shokrollahi1999}). In this case, in the fifth column of Table \ref{table2} we also reported the running time of the direct approach,
\emph{i.e.} using \eqref{basic-formula-L}, the third and fourth columns are respectively the running times
of the approaches \ref{Hasse-approach}) and \ref{Bernoulli-approach}).
For these values of $q$ it became clear that the computation time spent in performing the sums 
over $a$ was the longest one. This means that inserting an FFT-algorithm in the approaches  \ref{Hasse-approach}) and \ref{Bernoulli-approach}) is fundamental to further improve their performances.
We will say more on this in the next section.

\subsection{Computations  summing over $a$ via FFT (much faster, less decimal digits available).}
As we saw before,  as $q$ becomes large, the time spent in summing over $a$
dominates the overall computational cost.  So we implemented
the use of the FFT in  both the approaches \ref{Hasse-approach}) and \ref{Bernoulli-approach}),
by using  the {\tt fftw} \cite{FFTW} library in our   C programs.  
The performance of this part was extremely good in the sense that it was 
 a factor 1000 faster than the same one trivially performed. 

\subsection{Data for the scatter plots.}
\begin{figure} [h]
  \includegraphics[scale=0.95,angle=90]{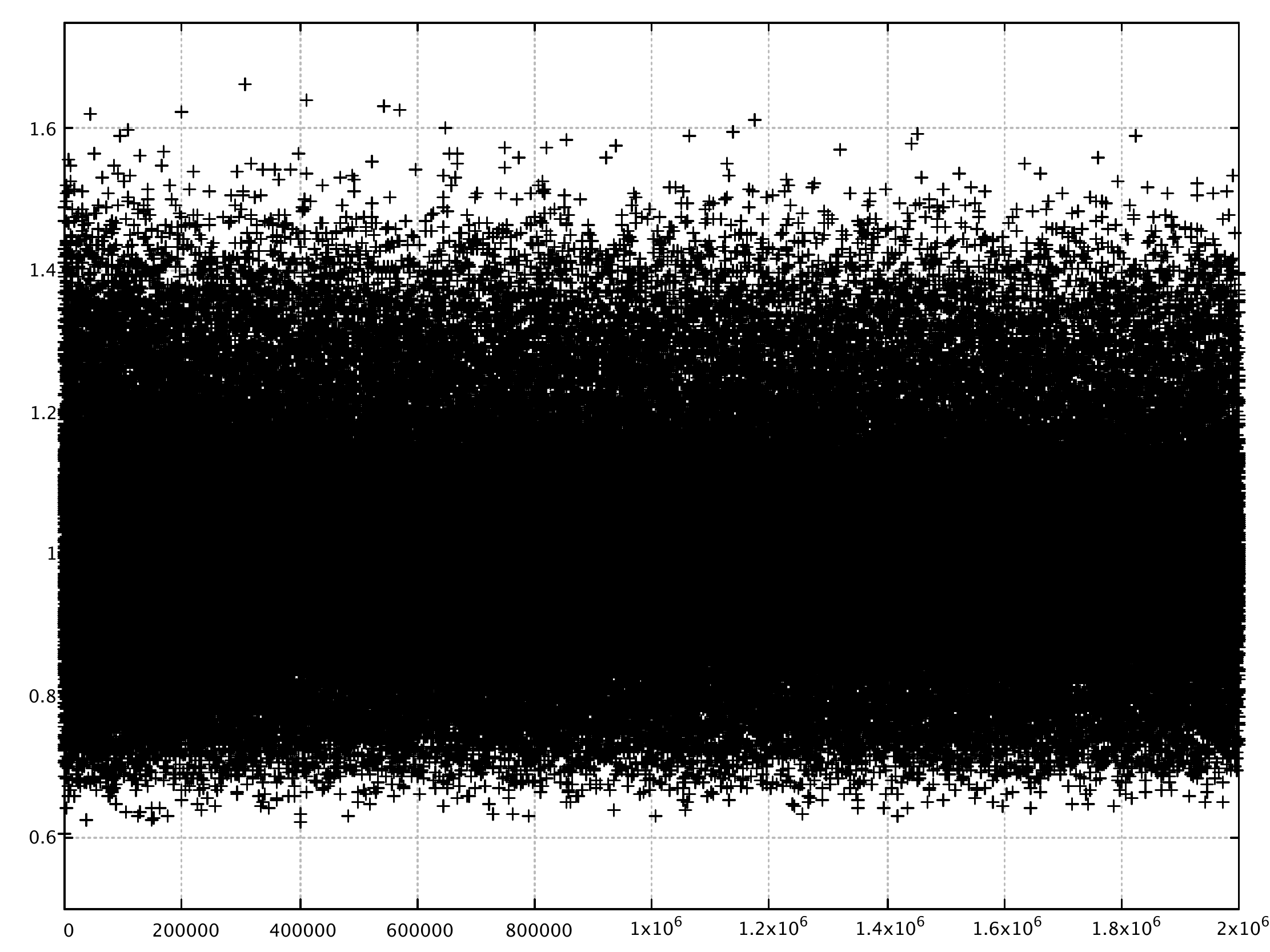}  
   \caption{{\small The values of $ r(q) $, $q$ prime, $3\le q\le  2\cdot 10^6$, plotted using
 GNUPLOT, v.5.2, patchlevel 7. [max $= 1.661436\dots$ 
 attained at $q= 305741$; min $=0.604599\dots$
 attained at $q= 3$;  number of $r(q)>1$: $74795$ ($50.22$\%); 
 number of $r(q)<1$: $74137$; ($49.78$\%) total number of  data $=148932$]. 
 }}
 \label{fig1}
 \end{figure}

We were able to compute the long double precision $r(q)$-values
 for every prime $3\le q \le  2\cdot 10^6$ and we provide here 
 the 
 scatter plot, see Figure \ref{fig1}, of  such values.
 The minimal value is $r(3) = 0.6045997880\dotsc$
 and the maximal one is   $r(305741) = 1.661436\dotsc$
The data were obtained  in about  three days of computation time on 
the Dell OptiPlex machine mentioned before. 

\subsection{Computations for larger $q$.}
If $bq+1$ is prime for many small $b,$ there is a good
chance that $r(q)$ will be large.
Promising, using this criterion, seemed
$q= 4178771$, $ 6766811$,  
$28227761$, $193894451$, $75743411$, $212634221$,
$251160191$, $405386081$, $538906601$, $964477901$ and also
$1139803271$, $1217434451$, $1806830951$, $2488788101$, 
$2830676081$, $2918643191.$ For these $q$
the $r(q)$'s were  evaluated using the  quadruple precision, see Table \ref{table3}.  
We remark that the quadruple precision  computation performances   are affected by a lack of hardware support
of the {\tt FLOAT128} type of the C programming language.  

In Table \ref{table4} we evaluated 
some further cases with potential large $r(q)$ value:
$q=4151292581$, $6406387241$, $7079770931$, 
$9854964401$ 
with the long double precision.
These computations were performed 
 on  an Intel(R) Xeon(R) CPU E5-2650 v3 @ 2.30GHz, with 160 GB of RAM
and running Ubuntu 16.04.
The computation for $q=9109334831$ was performed  
on the CAPRI  (``Calcolo ad Alte Prestazioni per la Ricerca e l'Innovazione'')
infrastructure of the University of Padova.

The PARI/Gp scripts and the C programs used and the computational results obtained
are available at the following web address:
\url{http://www.math.unipd.it/~languasc/rq-comput.html}. 
 
 \section{Further computation on the Euler-Kronecker constants} 
\begin{figure} 
 \includegraphics[scale=0.95,angle=90]{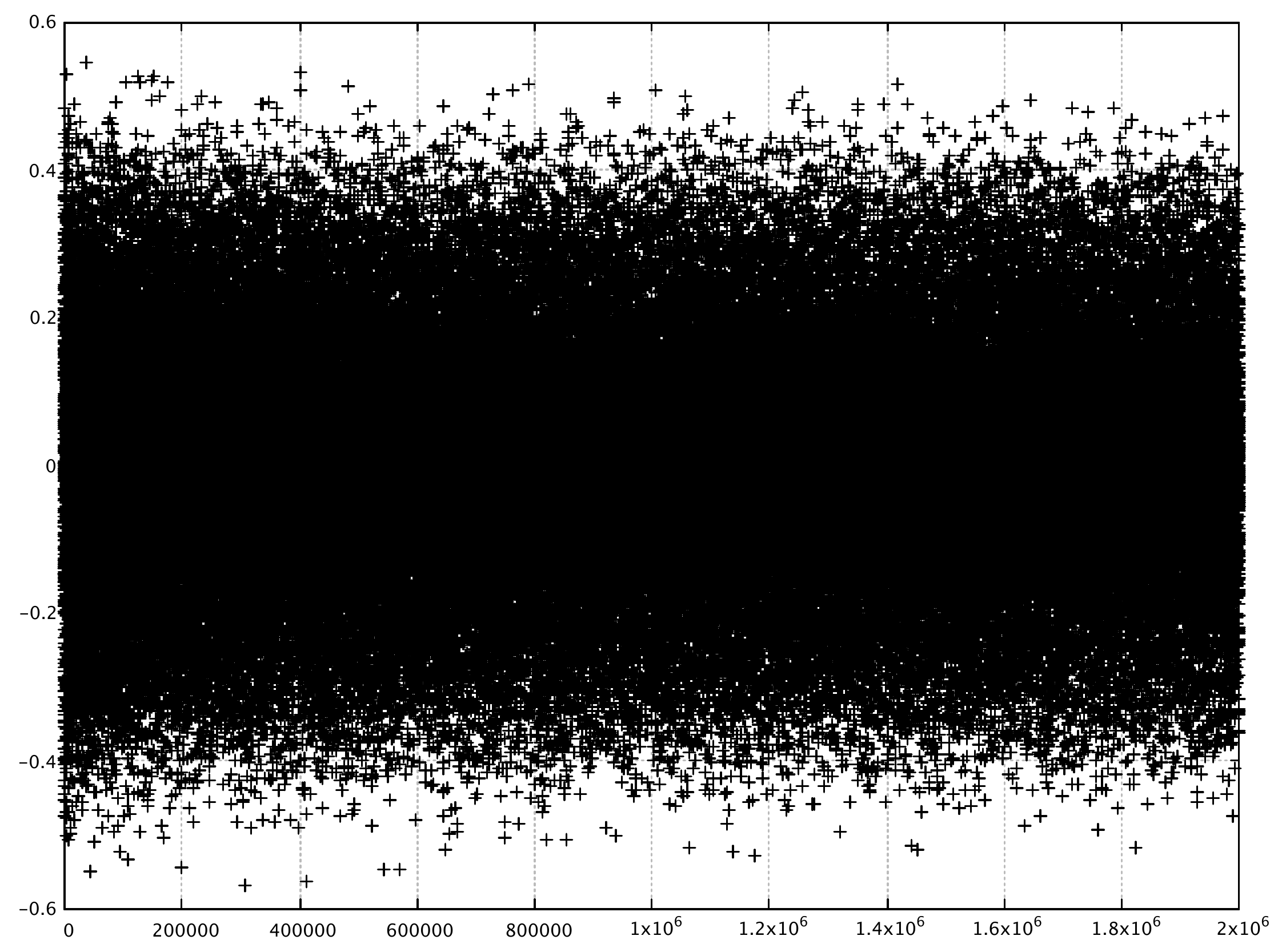}  
 \caption{\small{The values of $(\G_{q}-\G_{q}^+)/\log q$, $q$ prime, $3\le q\le  2\cdot 10^6$, plotted using
 GNUPLOT, v.5.2, patchlevel 7. [ 
max $=0.546473\dots$
 attained at $q= 37189$; min $=-0.569200\dots$
 attained at $q= 305741$; number of positive data $= 74190$  ($49.81$\%);
 number of negative data $=74742$ ($50.19$\%); total number of  data $=148932$]. 
}}
 \label{fig2}
 \end{figure}

 The Euler-Kronecker constant
for the prime cyclotomic field $\Q(\zeta_q)$
and for $\Q(\zeta_q+\zeta_q^{-1})$, the maximal real subfield
of $\Q(\zeta_q)$, see, \emph{e.g.}, \cite{Languasco2019},
are defined as
\begin{equation}
\label{EKq-EKq+def}
\G_{q}
: =
\gamma
+
\sum_{\chi \neq \chi_0} \frac{L'(1,\chi)}{L(1,\chi)},
\qquad
\textrm{respectively}
\qquad
\G_q^+ : = \gamma 
+
\sum_{\substack{\chi \neq \chi_0\\ \chi\, \textrm{even}}} 
\frac{L'(1,\chi)}{L(1,\chi)}.
\end{equation}
Using  formula  (22)  of \cite{Languasco2019} we have that
\begin{equation} 
\label{EKdiff-formula}
\G_{q}-\G_{q}^+ 
= 
\sum_{\chi\, \textrm{odd}} \frac{L'(1,\chi)}{L(1,\chi)}
=
\frac{q-1}{2}\Bigl(\gamma + \log(2\pi)\Bigr)
+
\sum_{\chi\, \textrm{odd}} 
 \frac{1}{B_{1,\overline{\chi}}}
 \sum_{a=1}^{q-1}  \overline{\chi}(a) \log\Bigl(\Gamma(\frac{a}{q})\Bigr),
\end{equation}
which can be easily implemented  since
the function $\log \Gamma$ is  available in the C programming language.
The
scatter plot of Figure \ref{fig2} represents the normalized values of
$\G_{q}-\G_{q}^+$, $q$ prime, $3\le q\le 2\cdot 10^6$ and 
creating it took about three days on
the Dell Optiplex machine mentioned before.
 
Since in \eqref{EKdiff-formula} we can use a decimation in frequency strategy, 
we also remark that letting $f = \log \Gamma$ into \eqref{DIF}  leads to simplify the form of $c_k
= e( - k / (q-1))
 \bigl(   \log \Gamma(a_k/q) -   \log \Gamma (a_{k+m}/q) \bigr)$, where $m=(q-1)/2$ and $k=0,\dotsc, m-1$,
in the following way\footnote{The minus sign here present in the \emph{twiddle factor} $e( - k / (q-1))$
comparing with the one in \eqref{DIF} depends on the fact that we are now summing over the conjugate Dirichlet
character $\overline{\chi}$ instead over $\chi$.}. Recalling $\langle g \rangle = \Z^*_q$, $a_k \equiv g^k \bmod q$ 
and   $g^m \equiv q-1 \bmod{q}$,
 we can write that 
\[ 
\log \Bigl(\Gamma \bigl(\frac{a_{k+m}}{q} \bigr)\Bigr)  = \log \Bigl(\Gamma\bigl(\frac{q-a_{k}}{q}\bigr)\Bigr)
=
\log \Bigl(\Gamma \bigl(1-\frac{a_{k}}{q}\bigr)\Bigr)
\] 
and hence, using the well-known \emph{reflection formula} $\Gamma(x) \Gamma(1-x)  = \pi / \sin(\pi x)$,  
we obtain
\begin{align*}
\log \Bigl(\Gamma \bigl(\frac{a_k}{q}\bigr)\Bigr)  - \log \Bigl(\Gamma \bigl(1-\frac{a_{k}}{q}\bigr)\Bigr)
&=
2 \log \Bigl(\Gamma \bigl(\frac{a_k}{q}\bigr)\Bigr)
+
\log\Bigl(\sin\bigl( \frac{\pi a_k}{q}\bigr)\Bigr) - \log \pi,
\end{align*}
for every $k=0,\dotsc, m-1$, $m=(q-1)/2$. Inserting the last relation in the definition of $c_k$  in \eqref{DIF}
we obtain the actual sequence we used for performing the computation previously mentioned.

\medskip 
\textbf{Acknowledgements}. 
Part of the calculations described here were performed
using the University of Padova Strategic Research Infrastructure Grant 2017:
``CAPRI: Calcolo ad Alte Prestazioni per la Ricerca e l'Innovazione''.
The first author, A.~Languasco, thanks Luca Righi (University of Padova) for his help in using   CAPRI.
The second and fourth authors, P.~Moree and A.~Sedunova, thank the Max Planck Institute for Mathematics
for providing excellent conditions while they were working on the project 
this paper is part of.
The third  author, S.~Saad Eddin, is supported by the Austrian Science Fund (FWF): Project F5505-N26 and Project F5507-N26, which are part of the special Research Program ``Quasi Monte Carlo Methods: Theory and Application''.


%
%
%


\newpage 
\newgeometry{left=0.75cm,right=0.5cm}
\begin{table}[htp]
\scalebox{0.6625}{
\begin{tabular}{|c|c|}
\hline
$q$  &  $r(q)$\\ \hline
$3$ & $0.60459978807807261686469275254\dotsc$\\ 
$5$ & $0.78956835208714868950675927999\dotsc$\\ 
$7$ & $0.95667518575084187547950733813\dotsc$\\ 
$11$ & $1.10916191287000575896982175317\dotsc$\\ 
$13$ & $1.07714905620985756748597815892\dotsc$\\ 
$17$ & $0.85539034568765268115905873936\dotsc$\\ 
$19$ & $0.70704004900384729070674621978\dotsc$\\ 
$23$ & $1.27303069939685502234405162961\dotsc$\\ 
$29$ & $1.19507225854723141702138692301\dotsc$\\ 
$31$ & $0.88988962107854407891985181571\dotsc$\\ 
$37$ & $0.89617354245182624263930105684\dotsc$\\ 
$41$ & $1.01095149281551337376703651618\dotsc$\\ 
$43$ & $1.00032807083987921579084335194\dotsc$\\ 
$47$ & $0.99510419475843763320461794597\dotsc$\\ 
$53$ & $1.00231549556080469808835403497\dotsc$\\ 
$59$ & $1.03111995957758588341749868917\dotsc$\\ 
$61$ & $0.91541689757636152038607840584\dotsc$\\ 
$67$ & $1.03230196304201968151553976333\dotsc$\\ 
$71$ & $0.94652474710362368092900546271\dotsc$\\ 
$73$ & $1.28217793230760538382246761185\dotsc$\\ 
$79$ & $0.84579459612002975504552940763\dotsc$\\ 
$83$ & $1.22326926548441461619507621390\dotsc$\\ 
$89$ & $1.28632147461922346234453694590\dotsc$\\ 
$97$ & $0.90467614287023765066781857933\dotsc$\\ 
$101$ & $1.11049958753586448051923888082\dotsc$\\ 
$103$ & $1.05565198833718743186163713482\dotsc$\\ 
$107$ & $0.99260767792672501309519612375\dotsc$\\ 
$109$ & $0.91554283885230186850667500246\dotsc$\\ 
$113$ & $1.16185573635061808057761114590\dotsc$\\ 
$127$ & $1.06269835499717635407980190888\dotsc$\\ 
$131$ & $1.27897699389762867270592988247\dotsc$\\ 
$137$ & $1.00188853650420792851571142833\dotsc$\\ 
$139$ & $0.87166115187392327886708542130\dotsc$\\ 
$149$ & $1.04886527642691194564791006449\dotsc$\\ 
$151$ & $1.09613526050530812035603232922\dotsc$\\ 
$157$ & $0.74304505329108896600523002862\dotsc$\\ 
$163$ & $0.95167392369442992883081838307\dotsc$\\ 
$167$ & $0.85404891714098835186838607451\dotsc$\\ 
$173$ & $1.25750311100604863256476652232\dotsc$\\ 
$179$ & $1.31898955218699008540672120548\dotsc$\\ 
$181$ & $1.01646725307901783240856438797\dotsc$\\ 
$191$ & $1.29850955347246763676155271715\dotsc$\\ 
$193$ & $1.17384956614280523683625176108\dotsc$\\ 
$197$ & $0.87142685805870225854275086741\dotsc$\\ 
$199$ & $0.79775765981803261703336410970\dotsc$\\ 
$211$ & $0.70965810384577007739153826881\dotsc$\\ 
$223$ & $0.90016736774009107389420074861\dotsc$\\ 
$227$ & $0.76298839763137122603762871170\dotsc$\\ 
$229$ & $0.72414574142010494620086404196\dotsc$\\ 
$233$ & $1.43102216731058063469583770264\dotsc$\\ 
$239$ & $1.18520259221018381028526578871\dotsc$\\ 
$241$ & $1.11908192699651325481120769078\dotsc$\\ 
$251$ & $1.18041694425392859170387583509\dotsc$\\ 
$257$ & $0.90559625735496576640913464538\dotsc$\\ 
$263$ & $0.93717078166852960654064932319\dotsc$\\ 
$269$ & $1.01052429941342866041104883014\dotsc$\\ 
\hline
\end{tabular}
}
\scalebox{0.6625}{
\begin{tabular}{|c|c|}
\hline
$q$  &  $r(q)$\\ \hline
$271$ & $0.84120880901441103034587178906\dotsc$\\ 
$277$ & $1.22287167700803659996325347046\dotsc$\\ 
$281$ & $1.09072312671446411507457756821\dotsc$\\ 
$283$ & $0.98730045924989351176735192970\dotsc$\\ 
$293$ & $1.28843023595237283191056798455\dotsc$\\ 
 $307$ & $0.91358725220199482220514916899\dotsc$ \\ 
 $311$ & $1.14589374542647302213444142687\dotsc$ \\ 
 $313$ & $0.93893317675819166180673984422\dotsc$ \\ 
 $317$ & $0.80671823188984812849457198577\dotsc$ \\ 
 $331$ & $0.81356274956051845902331649336\dotsc$ \\ 
 $337$ & $0.86111511521922591268832255791\dotsc$ \\ 
 $347$ & $1.08517941758105267446483313058\dotsc$ \\ 
 $349$ & $0.98395731344877010445591239132\dotsc$ \\ 
 $353$ & $0.88603505661744604503087815775\dotsc$ \\ 
 $359$ & $1.16002644446708254566916432735\dotsc$ \\ 
 $367$ & $0.90864101877936912063265319825\dotsc$ \\ 
 $373$ & $1.07507614420133257646267035535\dotsc$ \\ 
 $379$ & $0.72144618647138444694426995688\dotsc$ \\ 
 $383$ & $0.83243809267428713960470760380\dotsc$ \\ 
 $389$ & $0.84997782896854503971627563496\dotsc$ \\ 
 $397$ & $0.99757781120158579094253246616\dotsc$ \\ 
 $401$ & $1.13998328316447070631384278931\dotsc$ \\ 
 $409$ & $1.19919809743909540748744244798\dotsc$ \\ 
 $419$ & $1.18974458882376935926766971002\dotsc$ \\ 
 $421$ & $0.86457966530711741177342869535\dotsc$ \\ 
 $431$ & $1.13754261103593462461717085623\dotsc$ \\ 
 $433$ & $1.07176135182041771385450595205\dotsc$ \\ 
 $439$ & $0.68484134061729762055005895626\dotsc$ \\ 
 $443$ & $1.41089988430397986980906568345\dotsc$ \\ 
 $449$ & $0.90539643658614424895891547460\dotsc$ \\ 
 $457$ & $0.83734634190585621778636791343\dotsc$ \\ 
 $461$ & $1.03119557377397403645284724907\dotsc$ \\ 
 $463$ & $0.96134625111959841778686635231\dotsc$ \\ 
 $467$ & $0.89740454859192836870657083737\dotsc$ \\ 
 $479$ & $1.10506715780642069705910978939\dotsc$ \\ 
 $487$ & $1.13041022782656063139453697156\dotsc$ \\ 
 $491$ & $1.27221465691304968352754354985\dotsc$ \\ 
 $499$ & $0.82979024959465063669881382680\dotsc$ \\ 
  $503$ & $1.09956174719578329093362210466\dotsc$ \\ 
 $509$ & $1.39692082719612661320417410651\dotsc$ \\ 
 $521$ & $0.74488579181918272860910159248\dotsc$ \\ 
 $523$ & $0.99514847873992894203802693220\dotsc$ \\ 
 $541$ & $0.94472655782952981521345779529\dotsc$ \\ 
 $547$ & $0.73868505476195458996166613201\dotsc$ \\ 
 $557$ & $1.01800618130970440243478675140\dotsc$ \\ 
 $563$ & $0.92322125091337523644162001846\dotsc$ \\ 
 $569$ & $0.86644384514357385272705168284\dotsc$ \\ 
 $571$ & $0.99662480636851972762309151349\dotsc$ \\ 
 $577$ & $0.91370293804018510239277389204\dotsc$ \\ 
 $587$ & $0.81252459850672121660374173954\dotsc$ \\ 
 $593$ & $1.07734617489664930780759188442\dotsc$ \\ 
 $599$ & $0.96408773834723069779571268471\dotsc$ \\ 
  $601$ & $0.92827339751824097250854300550\dotsc$ \\ 
 $607$ & $0.83637312705251443247667799101\dotsc$ \\ 
 $613$ & $0.87703659303472148910355020294\dotsc$ \\ 
 $617$ & $0.84246084541946716141445378848\dotsc$ \\ 
\hline
\end{tabular}
}
\scalebox{0.6625}{
\begin{tabular}{|c|c|}
\hline
$q$  &  $r(q)$\\ \hline
 $619$ & $0.80463918636548231818097049238\dotsc$ \\ 
 $631$ & $1.13964698072442766479584447731\dotsc$ \\ 
 $641$ & $1.34299156432328475445263673245\dotsc$ \\ 
  $643$ & $1.01836205611360685304417553494\dotsc$ \\ 
 $647$ & $0.90233667317118875595490772209\dotsc$ \\ 
 $653$ & $1.27087727805772466468796098338\dotsc$ \\ 
 $659$ & $1.39106317898226550143998268526\dotsc$ \\ 
 $661$ & $0.83544430975232146569368385972\dotsc$ \\ 
 $673$ & $1.03660206982398637181187353211\dotsc$ \\ 
  $677$ & $0.92424013312497364401792044662\dotsc$ \\ 
 $683$ & $1.13528281402409476998254691134\dotsc$ \\ 
 $691$ & $0.76921427957454050696406411036\dotsc$ \\ 
 $701$ & $0.92089882867969861041626254382\dotsc$ \\ 
 $709$ & $1.05648934917801861606174800341\dotsc$ \\ 
 $719$ & $1.20306325855333927681117243729\dotsc$ \\ 
 $727$ & $0.99856921422780328631340639607\dotsc$ \\ 
 $733$ & $0.98014910177261986736078022621\dotsc$ \\ 
 $739$ & $1.10263546824053086630671245464\dotsc$ \\ 
 $743$ & $1.03495494096205775904091176831\dotsc$ \\ 
 $751$ & $1.01856200583585073878095848975\dotsc$ \\ 
 $757$ & $0.96706876118708598545541455447\dotsc$ \\ 
 $761$ & $1.46958285813141552491322656984\dotsc$ \\ 
 $769$ & $0.89892230367392111314972716474\dotsc$ \\ 
 $773$ & $1.06810947197031447130333305035\dotsc$ \\ 
 $787$ & $0.97178232843986336686451556474\dotsc$ \\ 
 $797$ & $1.03075130387360942943641986527\dotsc$ \\ 
 $809$ & $1.31970441406018712251949567645\dotsc$ \\ 
 $811$ & $0.80283817264815420707856818905\dotsc$ \\ 
 $821$ & $1.06528437036549643319352814656\dotsc$ \\ 
 $823$ & $0.96769318476182048656465705918\dotsc$ \\ 
 $827$ & $0.86555993675758442057691969953\dotsc$ \\ 
 $829$ & $0.82250033541615549748400919644\dotsc$ \\ 
 $839$ & $0.91871090540765761610044317666\dotsc$ \\ 
 $853$ & $1.08223582880253347548004283614\dotsc$ \\ 
 $857$ & $1.05075311490694694576396382027\dotsc$ \\ 
 $859$ & $0.88080094180568178476397675721\dotsc$ \\ 
 $863$ & $1.05694231206444764180240401289\dotsc$ \\ 
 $877$ & $0.72289398522705741218284637854\dotsc$ \\ 
 $881$ & $1.09738994199075350184435336371\dotsc$ \\ 
 $883$ & $1.13318227639393214982039012685\dotsc$ \\ 
 $887$ & $0.96917974196790823108417719936\dotsc$ \\ 
 $907$ & $0.90262558866311480478051623607\dotsc$ \\ 
 $911$ & $1.07798557536304873099351043701\dotsc$ \\ 
 $919$ & $1.04003346554199950901317303457\dotsc$ \\ 
 $929$ & $1.04414904452989167744813201725\dotsc$ \\ 
 $937$ & $0.90017934857750019784132262523\dotsc$ \\ 
 $941$ & $1.09400867179752235523397214847\dotsc$ \\ 
 $947$ & $1.22587448270510743026091490433\dotsc$ \\ 
 $953$ & $1.16083173031283885682226845601\dotsc$ \\ 
 $967$ & $0.72860004404668861481436825047\dotsc$ \\ 
 $971$ & $1.07939115916440046258710381604\dotsc$ \\ 
 $977$ & $0.83890885880371282354125472475\dotsc$ \\ 
 $983$ & $0.78867677202973854047246566763\dotsc$ \\ 
 $991$ & $0.90943936153505129760069639750\dotsc$ \\ 
 $997$ & $0.85575754491350654466545217865\dotsc$ \\ 
 \phantom{} & \phantom{}  \\ 
\hline
\end{tabular}
}
\caption{\label{table1}
Values of $r(q)$ for every odd prime up to $1000$ with 
a precision of $30$ digits; computed with PARI/Gp, v.~2.11.2,
with a trivial way of executing the sum over $a$.
Total computation time: 8s., 947ms.
using the first $\chi$-Bernoulli number. 
[s=seconds; msec=milliseconds]
}
\end{table} 
\restoregeometry 

\begin{table}[htp]
\begin{center}
\begin{tabular}{|c|c|r|r|r|}
\hline  
$q$  &  $r(q)$ \hfill &  time $\psi$-version &  time Bernoulli-version & time direct version \\ \hline
$1451$ & $1.489316072\dotsc$  & 292ms. & 299ms. & 2s. 083ms.\\
$2741$ & $1.498121015\dotsc$  & 1s. 019ms.& 1s. 061ms.& 5s. 231ms.\\
$3331$ & $0.642429297\dotsc$  & 1s. 496ms.& 1s. 566ms. & 7s. 236ms.\\
$4349$ & $1.518570512\dotsc$  & 2s. 536ms.& 2s. 657ms.& 10s. 572ms.\\
$4391$ & $1.507776410\dotsc$  & 2s. 586ms. & 2s. 710ms. & 11s. 276ms.\\
$5231$ & $1.556562248\dotsc$  & 3s. 361ms.& 3s. 385ms.& 14s. 730ms.\\
$6101$ & $1.511405291\dotsc$  & 4s. 964ms.& 5s.  086ms.& 17s. 490ms.\\
$6379$ & $0.673523026\dotsc$  & 5s. 436ms.& 5s.  699ms.& 19s. 859ms.\\
$7219$ & $0.658084090\dotsc$  & 6s. 978ms.& 7s.  293ms.& 23s. 966ms.\\
$8209$ & $0.672045039\dotsc$  & 8s. 950ms.& 9s.  416ms. & 27s. 322ms.\\
$9049$ & $0.667614244\dotsc$  & 10s. 870ms.&11s.  461ms. & 32s. 610ms.\\
$9689$ & $1.524371504\dotsc$  & 12s. 487ms.& 13s. 113ms.& 36s. 855ms.\\
\hline
\end{tabular}
\caption{\label{table2}
A few other values of $r(q)$   with 
a precision of $10$ digits; computed with PARI/Gp, v.~2.11.2,
with a trivial way of executing the sum over $a$.
[s=seconds; msec=milliseconds]
}
\end{center}
\end{table}

\begin{table}[htp]
\begin{center}
\begin{tabular}{|r|r|r|r|}  
\hline 
$q$ \hskip0.3truecm\mbox{}  &  $r(q)$  \hskip0.5truecm\mbox{}  & total time    &  total time     \\
    &      &  (long double)   &   (quadruple)   \\  \hline
$4178771$ & $1.611588128\dotsc$&2s. 471ms.&42s. 007ms.\\
$ 6766811$ & $1.709379041\dotsc$&4s. 200ms.&1m. 01s. 314ms.\\
$28227761$ & $1.528720351\dotsc$&16s. 267ms.& 4m. 30s. 804ms.\\ 
$75743411$ & $1.645759517\dotsc$&  1m. 05s. 078ms.& 19m. 04s. 250ms.\\ 
$193894451$ & $1.548501406\dotsc$& 2m. 18s. 739ms.& 38m. 27s. 838ms.\\ 
$212634221$ & $1.652149469\dotsc$&  2m.  44s. 018ms.&  48m. 24s. 277ms.\\ 
$251160191$ & $1.611898472\dotsc$& 3m. 02s. 783ms.& 54m.  20s.  875ms.\\ 
$405386081$&  $1.545118923\dotsc$& 4m. 31s. 634ms.& 58m.  49s.  163ms.\\ 
$538906601$ & $1.693680145\dotsc$& 6m. 44s. 584ms.& 114m. 59s. 330ms.\\ 
$964477901$ &   $1.612596619\dotsc$& 12m. 21s. 793ms.& 217m. 40s. 768ms.\\  
$1139803271$ & $1.398836497\dotsc$&  21m. 56s.  349ms.& 298m. 29s. 058ms.\\
$1217434451$ & $1.707310115\dotsc$& 15m. 00s. 374ms.& 277m. 58s. 390ms.\\  
$1806830951$ & $1.621464926\dotsc$ & 23m. 00s. 240ms.& 417m. 46s. 183ms.\\
$2488788101$ & $1.662760638\dotsc$ & 29m. 54s. 639ms.& 512m.   29s. 915ms.\\
$2830676081$ & $1.616923086\dotsc$ & 32m.   46s. 208ms.& 552m. 42s.  941ms.\\ 
$2918643191$ & $1.693092857\dotsc$& 38m. 40s. 423ms.& 697m. 10s. 305ms.\\  
 \hline
\end{tabular} 
\caption{\label{table3}
Few other values of $r(q)$; computed using the first $\chi$-Bernoulli number with 
PARI/Gp, v.~2.11.2. and fftw, v.~3.3.8., with  long double and quadruple  precisions. 
The sum over $a$ was performed using the FFT algorithm on the Xeon machine mentioned before.
[m=minutes; s=seconds; msec=milliseconds, n.a.=not available] 
}
\end{center}
\end{table}

\begin{table}[htp]
\begin{center}
\begin{tabular}{|r|r|r|}  
\hline 
$q$ \hskip0.3truecm\mbox{}  &  $r(q)$  \hskip0.5truecm\mbox{}  &  total time (long double)       \\ \hline 
$4151292581$ & $1.669735\dotsc$ &  63m. 00s.   \\
$6406387241$ & $1.625741\dotsc$ &  75m. 38s.  \\
$7079770931$& $1.688607\dotsc$ &   187m.  20s.  \\
$9109334831$ & $1.657855\dotsc$ &  120m. 28s. \\
$9854964401$ & $1.688033\dotsc$ &  177m. 14s. \\
 \hline
\end{tabular} 
\caption{\label{table4}
Few other values of $r(q)$; computed using the first $\chi$-Bernoulli number 
with PARI/Gp, v.~2.11.2. and fftw, v.~3.3.8., with  long double  precision. 
The sum over $a$ was performed using the FFT algorithm on the Xeon machine mentioned before
with the only exception of $q=9109334831$ for which we used the  CAPRI infrastructure of the
University of Padova. [m=minutes; s=seconds] 
}
\end{center}
\end{table}  

\phantom{1223}

\end{document}